\documentclass[a4paper,12pt,reqno]{amsart}

\usepackage{color}
\usepackage{amsmath}
\usepackage{amssymb}
\usepackage{amsfonts}
\usepackage{graphicx}
\usepackage{mathtools}
\usepackage[colorlinks]{hyperref}
\renewcommand\eqref[1]{(\ref{#1})} 
\usepackage{booktabs}
\usepackage{lscape}

\usepackage{mathrsfs}
\usepackage{comment}

\graphicspath{ {images/} }
\newcommand{\Q}{\mathbb{Q}}
\newcommand{\R}{\mathbb{R}}

\title[Heat and wave type equations]{$L^p-L^q$ estimates for non-local heat and wave type equations on locally compact groups}
\author[S. G\'omez Cobos]{Santiago G\'omez Cobos}
\address{
	Santiago G\'omez Cobos:
	\endgraf
	Department of Mathematics: Analysis, Logic and Discrete Mathematics
	\endgraf
	Ghent University, Krijgslaan 281, Building S8, B 9000 Ghent
	\endgraf
	Belgium
	\endgraf
	{\it E-mail address} {\rm davidsantiago.gomezcobos@ugent.be}}

\author[J. E. Restrepo]{Joel E. Restrepo}
\address{
	Joel E. Restrepo:
	\endgraf
	Department of Mathematics: Analysis, Logic and Discrete Mathematics
	\endgraf
	Ghent University, Krijgslaan 281, Building S8, B 9000 Ghent
	\endgraf
	Belgium
	\endgraf
	{\it E-mail address} {\rm joel.restrepo@ugent.be; cocojoel89@yahoo.es}}

\author[M. Ruzhansky]{Michael Ruzhansky}
\address{
	Michael Ruzhansky:
	\endgraf
	Department of Mathematics: Analysis, Logic and Discrete Mathematics
	\endgraf
	Ghent University, Krijgslaan 281, Building S8, B 9000 Ghent
	\endgraf
	Belgium
	\endgraf
	and
	\endgraf
    School of Mathematical Sciences
    \endgraf
    Queen Mary University of London
    \endgraf
    United Kingdom
    \endgraf
	{\it E-mail address} {\rm michael.ruzhansky@ugent.be}}


\thanks{The paper has been accepted for publication in the journal ``Comptes Rendus Math\'{e}matique". }

\newtheoremstyle{theorem}
{10pt}          
{10pt}  
{\sl}  
{\parindent}     
{\bf}  
{. }    
{ }    
{}     
\theoremstyle{theorem}

\numberwithin{equation}{section}
\theoremstyle{plain}
\newtheorem{thm}{Theorem}[section]

\theoremstyle{definition}
\newtheorem{defn}[thm]{Definition}

\newtheoremstyle{defi}
{10pt}          
{10pt}  
{\rm}  
{\parindent}     
{\bf}  
{. }    
{ }    
{}     
\theoremstyle{defi}

\begin{document}

 	\begin{abstract}

We prove the $L^p-L^q$ $(1<p\leqslant 2\leqslant q<+\infty)$ norm estimates for the solutions of heat and wave type equations on a locally compact separable unimodular group $G$ by using an integro-differential operator in time and any positive left invariant operator (maybe unbounded) on $G$. We complement our studies by  giving asymptotic time estimates for the solutions, which in some cases are sharp. 

\end{abstract}
	\maketitle
	\tableofcontents

\section{Introduction}

As an application of some spectral multipliers results, it was shown recently in \cite{RR2020} that the  $L^p-L^q$ norm estimates for the solution of the $\mathscr{L}$-heat equation 
\begin{align*}
\begin{split}
\partial_{t}w(t,x)+\mathscr{L}w(t,x)&=0, \quad t>0,\,\, x\in G, \\
w(t,x)|_{_{_{t=0}}}&=w_0(x),
\end{split}
\end{align*}
on a locally compact separable unimodular group $G$ can be reduced to the time asymptotics of its propagator in the noncommutative Lorentz space norm \cite{[51]}, which involves calculating the trace of the spectral projections of the operator $\mathscr{L}$. The considered operator can be any positive linear left invariant operator acting on $G$ with the possibility (generality) of having either continuous or discrete spectrum. From convention of \cite{RR2020}, a positive operator means nonnegative and self-adjoint. Unfortunately, the results did not cover the $L^p-L^{q}$ estimates for the case of $\mathscr{L}$-wave equations. The propagators can not be estimated by the noncommutative Lorentz space norm. This is still an open question. In this paper, we fill this gap if we consider a non-local integro-differential operator in time of order $\beta$ $(0<\beta<2)$, which allows us to consider $\mathscr{L}$-equations which interpolate between heat and wave types.

The type of groups that we can consider is very large. For instance, compact, semi simple, exponential, nilpotent, some solvable ones, real algebraic, and many more.

Our results are new in this context and open a door to studying different estimates of other types of equations which, in principle, can not be treated in the classical case but perhaps could be carried out in the set up of non-local operators.   

Thus the main results of the paper concern $L^p-L^q$ $(1<p\leqslant 2\leqslant q<+\infty)$ norm estimates for $\mathscr{L}$-heat-wave type equations. Full details are given in Section \ref{heat-locally}. We then give a number of examples by choosing concrete operators and groups. Moreover, in some cases, we are able to claim the sharpness of the time-decay rate. 

\section{Von Neumann algebras and noncommutative Lorentz spaces}\label{preli}

Let $\mathfrak{L}(\mathcal{H})$ be the set of linear operators defined on a Hilbert space $\mathcal{H}$. In this context, the concept of $\tau$-measurability on a von Neumann algebra $M$ (see e.g. \cite{von,[46],[47],von2}) and the spectral projections give us the possibility to approximate unbounded operators by bounded ones. 

Let us now recall some important definitions which will be used frequently in the development of this paper.   

\begin{defn}
Let $M\subset\mathfrak{L}(\mathcal{H})$ be a semifinite von Neumann algebra acting over the Hilbert space $\mathcal{H}$ with a trace $\tau$. A linear operator $L$ (maybe unbounded) is said to be affiliated with $M$, if it commutes with the elements of the commutant $M^!$ of $M$, i.e. $LV=VL$ for all $V\in M^!.$
\end{defn}  

\begin{defn}
A closeable operator $L$ (maybe unbounded) is called  $\tau$-measureable if for each $\epsilon>0$ there exists a projection $P$ in $M$ such that $P(\mathcal{H})\subset D(L)$ and $\tau(I-P)\leqslant \epsilon$, where $D(L)$ is the domain of $L$ in $\mathcal{H}$ and $M\subset\mathfrak{L}(\mathcal{H})$. We denote by $S(M)$ the set of all $\tau$-measurable operators.   
\end{defn}

\begin{defn}
Let $L\in S(M)$, and let $L = U|L|$ be its polar decomposition. We define the distribution function by $d_\gamma(L):=\tau\big(E_{(\gamma,+\infty)}(|L|)\big)$ for $\gamma\geqslant0$, where $E_{(\gamma,+\infty)}(|L|)$ is the spectral projection of $L$ over the interval $(\gamma,+\infty).$ Also, for any $t>0$, we define the generalized $t-$th singular numbers as 
\[
\mu_t(L):=\inf\{\gamma\geqslant0:\,d_\gamma(L)\leqslant t\}.
\]
\end{defn}
More discussions (and properties) of the distribution function and generalized singular numbers can be found e.g. in \cite{pacific}. 

We conclude this section by recalling the noncommutative Lorentz spaces associated with a semifinite von Neumann algebra $M$ \cite{[51]}.

\begin{defn}
We denote by $L^{p,q}(M)$ $(1\leqslant p<+\infty,\,1\leqslant q<+\infty)$ the set of all operators $L\in S(M)$ such that
\begin{align*}
\|L\|_{L^{p,q}(M)}=\left(\int_0^{+\infty}\big(t^{1/p}\mu_t(L)\big)^q \frac{{\rm d}t}{t}\right)^{1/q}<+\infty.    
\end{align*}
Therefore the $L^{p}$-spaces on $M$ are given by
\[
\|L\|_{L^{p}(M)}:=\|L\|_{L^{p,p}(M)}=\left(\int_0^{+\infty}\mu_t(L)^p{\rm d}t\right)^{1/p}.
\]
For the case $q=\infty$, $L^{p,\infty}(M)$ is the set of all operators $L\in S(M)$ with  
\[
\|L\|_{L^{p,\infty}(M)}=\sup_{t>0}t^{1/p}\mu_t(L).
\]
\end{defn}   

\section{Heat and wave type equations}\label{heat-locally}

In this section we establish for non-local heat and wave type equations the $L^p-L^q$ norm estimates for $1<p\leqslant 2\leqslant q<+\infty$. Indeed, we will see that the norm estimates can be reduced to the time asymptotics of the propagator in the noncommutative Lorentz space norm. Also, we show that the latter norm mainly involves calculating the trace of the spectral projections of the operator $\mathscr{L}$.

In the presentation of the propagators of this section, frequently the two-parametric Mittag-Leffler function will appear:
\[
E_{\alpha,\delta}(z)=\sum_{k=0}^{+\infty} \frac{z^k}{\Gamma(\alpha k+\delta)},\quad z,\delta\in\mathbb{C},\quad \Re(\alpha)>0,
\]
which is absolutely and locally uniformly convergent for the given parameters (\cite{mittag}).

\medskip For the existence, uniqueness, representation and other properties of solutions for the considered equations below, we refer to the theory on evolutionary integral equations in \cite[Chapter I]{pruss} or, alternatively, one can use the Borel functional calculus \cite{BorelFunctional} to construct the propagators.  More concrete and similar studies on abstract spaces can be found in \cite[Chapter 3]{thesis} (see also \cite{thesis2001}). For the estimation of the $L^p-L^q$ norms, we usually combine \cite[Theorem 5.1]{RR2020} with the alternative form provided by \cite[Theorem 6.1]{RR2020} to calculate the noncommuative Lorentzian norm.

\subsection{\texorpdfstring{$\mathscr{L}$}{L}-heat type equation} We consider the following heat type equation: 
\begin{equation}\label{heatlocally}
\begin{split}
\underbrace{\int_0^t \frac{(t-s)^{-\beta}}{\Gamma(1-\beta)}\partial_{s}w(s,x)\,\mathrm{d}s}_{^{C}\partial_{t}^{\beta}w(t,x)}+\mathscr{L}w(t,x)&=0, \quad t>0,\,\, x\in G, \\
w(t,x)|_{_{_{t=0}}}&=w_0(x),
\end{split}
\end{equation}
where $^{C}\partial_{t}^{\beta}$ is the so-called Dzhrbashyan-Caputo fractional derivative, $0<\beta<1$, $G$ is a locally compact separable unimodular group and $\mathscr{L}$ is any positive linear left invariant operator on $G$ (maybe unbounded). Our next result also covers the case $\beta=1$, i.e. the classical $\mathscr{L}$-heat equation, but it is already known \cite[Section 7]{RR2020}. 

\begin{thm}\label{heatlocally-thm}
Let $G$ be a locally compact separable unimodular group, $0<\beta<1$ and $1<p\leqslant 2\leqslant q<+\infty$. Let $\mathscr{L}$ be any positive left invariant operator on $G$ (maybe unbounded) such that
\[
\sup_{t>0}\sup_{s>0}\big[\tau\big(E_{(0,s)}(\mathscr{L})\big)\big]^{\frac{1}{p}-\frac{1}{q}}E_\beta(-t^{\beta}s)<+\infty. 
\]
If $w_0\in L^p(G)$ then there exists a unique solution $w\in\mathcal{C}\big([0,+\infty) ;L^q(G)\big)$ for the $\mathscr{L}$-heat type equation \eqref{heatlocally} represented by
\[
w(t,x)=E_\beta(-t^{\beta}\mathscr{L})w_0(x),\quad t>0,\,\,x\in G.
\]
In the particular case that  
\begin{equation}\label{tracecondition}
\tau\big(E_{(0,s)}(\mathscr{L})\big)\lesssim s^{\lambda},\quad s\to+\infty,\quad\text{for some}\quad \lambda>0,
\end{equation}
we get the following time decay rate for the solution of equation \eqref{heatlocally} for all $t>0,$
\[
\|w(t,\cdot)\|_{L^q(G)}\leqslant C_{\beta  ,\lambda,p,q}t^{-\beta\lambda\left(\frac{1}{p}-\frac{1}{q}\right)}\|w_0\|_{L^p(G)},\quad \frac{1}{\lambda}\geqslant\frac{1}{p}-\frac{1}{q}.    
\]
\end{thm}

\subsection{\texorpdfstring{$\mathscr{L}$}{L}-wave type equation}\label{wave-locally}

Here we investigate the solution of the following equation, which in a sense interpolates between wave (without being wave, $\beta<2$) and heat types:  
\begin{equation}\label{locallywave}
\begin{split}
^{C}\partial_{t}^{\beta}w(t,x)+\mathscr{L}w(t,x)&=0, \quad t>0,\,\, x\in G, \\
w(t,x)|_{_{_{t=0}}}&=w_0(x), \\
\partial_t w(t,x)|_{_{_{t=0}}}&=w_1(x).
\end{split}
\end{equation}
where $^{C}\partial_{t}^{\beta}$ is the Dzhrbashyan-Caputo fractional derivative,  $\mathscr{L}$ is any positive linear left invariant operator on $G$ (maybe unbounded) and $1<\beta<2$. 

\medskip Before giving the main result of this subsection, first we mention a technical and important theorem whose proof will appear somewhere else. Remember that for $L\in S(M)$ we denote $L=U|L|$ its polar decomposition.

\begin{thm}\label{additional}
Let $L$ be a closed (maybe unbounded) operator affiliated with a semifinite von Neuman algebra $M$. Let $\phi$ be a Borel measurable function on $[0,+\infty)$. Suppose also that $\psi$ is a monotonically decreasing continuous function on $[0,+\infty)$ such that $\psi(0)=1$, $\lim_{v\to+\infty}\psi(v)=0$ and $\textcolor{red}{|}\phi(v)\textcolor{red}{|}\leqslant \psi(v)$ for all $v\in[0,+\infty).$ Then for every $1\leqslant r<\infty$ we have the inequality 
\[
\|\phi(|L|)\|_{L^{r, \infty}(M)}\leqslant \sup_{v>0}\psi(v)\big[\tau(E_{(0,v)}(|L|))\big]^{\frac{1}{r}}.
\]
\end{thm}
Utilizing Theorem \ref{additional} we will be able to prove and state the main result on wave type equations. Now we need to recall an operator which is involved in the representation of the wave type propagator. Therefore, we introduce $\prescript{RL}{0}I^{\beta}u(t)=\frac1{\Gamma(\beta)}\int_0^t (t-s)^{\beta-1}u(s)\,\mathrm{d}s$, the Riemann-Liouville fractional integral of order $\beta>0,$ which is well-defined for functions on $L^1(0,T).$

\begin{thm}\label{locally-wave-thm}
Let $G$ be a locally compact separable unimodular group, $1<\beta<2$ and $1<p\leqslant 2\leqslant q<+\infty$. Let $\mathscr{L}$ be any positive left invariant operator on $G$ (maybe unbounded) satisfying the condition  
\[
\sup_{t>0}\sup_{s>0}\frac{\big[\tau\big(E_{(0,s)}(\mathscr{L})\big)\big]^{\frac{1}{p}-\frac{1}{q}}}{1+t^\beta s}<+\infty.
\]
If $w_0,w_1\in L^p(G)$ then there exists a unique solution $w\in\mathcal{C}\big([0,+\infty);L^q(G)\big)$ for the $\mathscr{L}$-wave type equation \eqref{locallywave} given by
\[
w(t,x)=E_\beta(-t^{\beta}\mathscr{L})w_0(x)+\prescript{RL}{0}I^{1}_t E_{\beta}(-t^{\beta}\mathscr{L})w_1(x),\quad t>0,\,\,x\in G.
\]
In particular, if the condition \eqref{tracecondition} holds then for any $1<p\leqslant 2\leqslant q<+\infty$ such that $\frac{1}{\lambda}\geqslant\frac{1}{p}-\frac{1}{q}$ we get the following time decay rate for the solution of equation \eqref{locallywave}: 
\[
\|w(t,\cdot)\|_{L^q(G)}\leqslant  C_{\beta,\lambda,p,q}t^{-\beta\lambda\left(\frac{1}{p}-\frac{1}{q}\right)} \big(\|w_0\|_{L^p(G)}+t\|w_1\|_{L^p(G)}\big).    
\]
\end{thm}
In Theorems \ref{heatlocally-thm}, \ref{locally-wave-thm}, the time decay rate for the solutions of equations \eqref{heatlocally} and \eqref{locallywave}, is predetermined by the condition \eqref{tracecondition}. Hence, let us mention briefly several examples of operators (on different groups) such that the trace of the spectral projections behave like $s^{\lambda}$ as $s\to+\infty$. In the Euclidean space $\R^n$, consider the Laplacian $\Delta_{\R^n}$. Here, we have such a behaviour \cite[Example 7.3]{RR2020} with $\lambda=n/2.$ Moreover it is given in more generality since it allows one to consider not just conjugate exponents. Here we can recover the sharp estimate (time-decay) given in \cite[Theorem 3.3, item (i)]{uno2} whenever $\frac{2}{n}>\frac{1}{p}-\frac{1}{q}$. The sub-Laplacian $\Delta_{sub}$ on a compact Lie group \cite{[35]} with $\lambda=Q/2$, where $Q$ is the Hausdorff dimension of $G.$ The positive sub-Laplacian on the Heisenberg group $\mathbb{H}^n$ \cite[Formula (7.17)]{RR2020} with $\lambda=n+1.$ A positive Rockland operator $\mathcal{R}$ of homogeneous order $\nu$ on a graded Lie group \cite[Theorem 8.2]{david} with $\lambda=Q/\nu$, where $Q$ is the homogeneous dimension of $G$. By using the latter estimate we can get sharp time decay in Theorems \ref{heatlocally-thm}, \ref{locally-wave-thm}. For operators of non-Rockland-type on the Engel group $\mathfrak{B}_4$ $(\lambda=3)$ and the Cartan group $\mathfrak{B}_5$ $(\lambda=9/2)$ \cite[Examples 2.2 and 3.2]{Marianna}. An $m$-th order weighted subcoercive positive operator on a connected unimodular Lie group \cite[Proposition 0.3]{david2} with $\lambda=Q^*/m$, where $Q_*$ is the local dimension of $G$ relative to the chosen weighted structure on its Lie algebra.

\section{Acknowledgements}
The authors were supported by the FWO Odysseus 1 grant G.0H94.18N: Analysis and Partial Differential Equations, the Methusalem programme of the Ghent University Special Research Fund (BOF) (Grant number 01M01021). MR is also supported by EPSRC grant EP/R003025/2 and FWO Senior Research Grant G011522N.

\end{document}